\documentclass[aps,preprint]{elsarticle}%
\usepackage{amsfonts}
\usepackage{amsmath}
\usepackage{amssymb}
\usepackage{graphicx}%
\providecommand{\U}[1]{\protect\rule{.1in}{.1in}}

\begin{document}

\title{Convergence analysis and parity conservation of a new form of a quadratic explicit spline}
\author{A.J. Ferrari$^{a,b}$, L.P. Lara$^{a,b}$, E.A. Santillan Marcus$^b$}

\address{$^a$ Consejo Nacional de Investigaciones Cient\'ificas y T\'ecnicas, 27 de febrero 210 bis, S2000, Rosario, Argentina. \\ $b$ Departamento de Matem\'atica, Facultad de Ciencias Exactas, Ingenier\'ia y Agrimensura, Universidad Nacional de Rosario, Av. Pellegrini 250, S2000, Rosario, Argentina.}

\maketitle

\textbf{Abstract.} In this study, a new form of quadratic spline is obtained, where the coefficients are determined explicitly by variational methods. Convergence is studied and parity conservation is demonstrated. Finally, the method is applied to solve integral equations. \\


\section{INTRODUCTION}

In Mathematics, Physics and Engineering, among other disciplines, there is a great need to adjust a discrete set of data or to approximate functions. In general, it is desired to know the values ​​at intermediate points, which can be solved through interpolation polynomials. In the practice, high-order polynomials can introduce significant errors due to various factors. Generally, these polynomials present fluctuations that are not present in the function to interpolate. For this reason, in this work the piecewise interpolation will be considered, which is particularly useful when the data to adjust have a smooth behavior alternated with strong changes. Focus will be on quadratic piecewise interpolation $S$ of continuous real functions. In \cite{kincaid91, henrici64, press86} methods and algorithms that depend on a determined criterion in the determination of one of the coefficients of $S$ are presented. Usually, all these methods require the resolution of algebraic systems or recursive equations. In our case, a
variational alternative that minimizes the fluctuations of the interpolator polynomial $S$ is presented. Its coefficients are determined explicitly through elementary functions. This work is organized as follows: the interpolator is presented in section 2; the convergence is studied in section 3; it is shown that $S$ maintains the parity of the function to be interpolated in section 4; the results are applied for Fredholm linear integral equations in section 5; similarly for Volterra linear integral equations in section 6; the numerical results are shown in section 7; and finally the conclusions are presented in section 8.

\section{THE QUADRATIC SPLINE}

Consider an interval $[a,b]$ in which we select $n+1$ equidistant nodes $x_{0}=a<x_1<\ldots<x_n=b$, where $x_k-x_{k-1}=h$,
$k=1,\ldots,n$ . Let $y:[a,b]\rightarrow\mathbb{C}$ be a continuous function, 
and consider $y_k=y(x_{k})$, $k=0,\ldots,n$. Let $I_k:=[x_{k-1},x_k]$, $k=1,\ldots,n.$ It is desired to determine the
quadratic piecewise interpolator $S(x)$ such that interpolates $y(x)$ in $x_k$, $k=0,\ldots,n$ and such that
$S^\prime(x)$ is continuous in the nodes $x_k$, $k=1,\ldots,n-1$.

Let $S(x)$ be the function in $[x_0,x_n]$
defined through the polynomials $P_k(x)$ so that:
\begin{equation} \label{1}
S(x)=P_k(x), \; x \in I_k, \; k=1,\ldots,n.
\end{equation}

For being $S$ interpolant, $P_k(x)$ must verify:
\begin{align}
P_k(x_{k-1})  &  = y_{k-1}, \; k=1,\ldots,n, \nonumber \\
P_k(x_k)  &  = y_k, \; k=1,\ldots,n.\nonumber
\end{align}

The continuity of $S^{\prime}(x)$ imposes that:
\begin{equation} \label{2bis}
P_k^\prime(x_k)=P_{k+1}^\prime(x_k), \; k=1,\ldots,n-1.
\end{equation}

Lagrange polynomials $p_{k}(x)$ over $I_{k}$ are built as:
\begin{equation} \label{3bis}
p_k(x)=\frac{x-x_{k-1}}{h} \, y_{k}-\frac{x-x_{k}}{h} \, y_{k-1}, \; k=1,\ldots,n,
\end{equation}

\noindent which satisfy $p_k(x_{k-1})=y_{k-1}$, $p_k(x_k)=y_k$, $k=1,\ldots,n.$ Then:
\begin{equation} \label{4}
P_k(x)=p_k(x)+a_k(x-x_{k-1})(x-x_k), \; k=1,\ldots,n.
\end{equation}

In \cite{foucher09, rana89} the coefficients $a_k$ are obtained by linear algebraic systems. Here an explicit formula will be found.
From (\ref{2bis}) and through simple algebraic operations is obtained:%
\begin{align} \label{4bis}
a_{k+1}  &  =\Delta_k-a_k, \; k=1,\ldots,n-1, \\
\Delta_k  &  =\frac{y_{k-1}-2y_k+y_{k+1}}{h^2} \nonumber,
\end{align}

\noindent from where $a_k=a_k(a_{1})$, linear in $a_1$. The coefficient $a_1$ is determined from an additional condition.

From (\ref{4bis}):
\begin{equation} \label{444new}
a_k=(-1)^{k+1}a_1+r_k, \; k \geq 2,
\end{equation}
\begin{equation} \label{4444}
r_k=(-1)^{k+1}\sum\limits_{j=1}^{k-1}(-1)^j\Delta_j, \; k \geq 2.
\end{equation}

Thus:
\begin{equation} \label{4445}
a_k=(-1)^{k+1}\left(a_1+\sum\limits_{j=1}^{k-1}(-1)^{j}\Delta_j\right), \; k \geq 2.
\end{equation}

In this way an explicit expression is available in terms of
elementary functions, to calculate the coefficients of the polynomial $S$. In addition it is simple to show that:
\begin{align}
r_{k+1}  &  =\Delta_k-r_{k}, \; k=1,\ldots,n-1, \nonumber \\
r_1  &  = 0. \nonumber
\end{align}

Next, $a_{1}$ is determined so that the sum of the quadratic errors between $P_k(x)$ and $p_k(x)$ in each $I_k$ be minimum. Let $E$ be defined as $E=\sum_{k=1}^n E_k$, where:
\begin{equation} \nonumber
E_k=\int\limits_{x_{k-1}}^{x_k}[P_k(x)-p_k(x)]^2\,dx.
\end{equation}

Taking into account (\ref{3bis},\ref{4}), we have:
\begin{equation} \nonumber
E=E(a_1)=\sum\limits_{k=1}^{n}a_{k}^2\int\limits_{x_{k-1}}^{x_k}(x-x_k)^2(x-x_{k-1})^{2}\,dx,
\end{equation}

\noindent from where:
\begin{equation} \nonumber
E(a_1)=\frac{h^5}{30}\sum\limits_{k=1}^n a_k^2.
\end{equation}

Then, taking into account (\ref{444new}), results $\frac{\partial}
{\partial a_1}a_k=(-1)^{k+1}$ from where $\frac{\partial^2 E(a_1)}{\partial a_1^2}=\frac{1}{15} n h^5 > 0$. Therefore, the solution $a_1$ that minimizes $E$ is such that $\frac{\partial E(a_1)}{\partial a_1}=0$. This leads to:
\begin{equation} \nonumber
a_1=\frac{1}{n}\sum\limits_{k=1}^n (-1)^k r_k.
\end{equation}

From (\ref{4444}) arises:
\begin{equation} \label{666}
a_1=-\frac{1}{n}\sum\limits_{j=1}^{n-1}(n-j)(-1)^j\Delta_{j}.
\end{equation}

Finally combining (\ref{4445}) and (\ref{666}):
\begin{equation} \label{ak}
a_k=(-1)^{k+1}\sum\limits_{j=1}^{n-1}\left(\frac{j}{n}+s_j-1\right)(-1)^j\Delta_j, \; k \geq 2,
\end{equation}

\noindent where $s_j=1$ if $j \leq k-1$ and $s_j=0$ if $j>k-1$.

Taking into account the definition of $\Delta_j$:
\begin{equation} \label{1313}
a_k=\sum_{j=0}^n c_{k,j} \, y_{j},
\end{equation}
where $c_{k,j}$ is given by:

\begin{equation} \nonumber
c_{k,j}=
\left\{\begin{array}{ll}
\frac{(-1)^k}{h^2}\beta_1 & \text{if} \; j=0 \\
\frac{(-1)^{k+1}}{h^2}(2\beta_1+\beta_2) & \text{if} \; j=1 \\
\frac{(-1)^{k+j}}{h^2}(\beta_{j-1}+2\beta_j+\beta_{j+1}) & \text{if} \; 1<j<n-1 \\
\frac{(-1)^{k+n-1}}{h^2}(2\beta_{n-2}+\beta_{n-1}) & \text{if} \; j=n-1 \\
\frac{(-1)^{k+n}}{h^2}\beta_{n-1} & \text{if} \; j=n \\
\end{array}\right.
\end{equation}

\noindent where $\beta_j=\frac{j}{n}$ if $j \leq k-1$ and $\beta_j=\frac{j}{n}-1$ if $j > k-1$.

(\ref{1313}) can be rewritten in matrix form as $A=C \cdot Y$ being the matrixes $\mathit{Y}_{n+1,1}=\{y_j\}_{j=0,\ldots,n}$ and
$\mathit{C}_{n,n+1}=\{c_{k,j}\}_{k=1,\ldots,n}^{j=0,\ldots,n}$. Note that the latter only depends
on $x_{0}$, $h$, and $n$. (\ref{3bis}) and (\ref{4}) define the matrixes $\mathit{p}_{n,1}(x)=\{p_k(x)\}_{k=1,\ldots,n}$, $\mathit{P}_{n,1}(x)=\{P_{k}(x)\}_{k=1,\ldots,n}$. Then the spline $S$ can be written in the following matrix form:
\begin{equation} \nonumber
\mathit{P}(x)=\mathit{p}(x)+\mathit{X(x) \cdot A} = \mathit{p}(x)+\mathit{X(x) \cdot C \cdot Y}.
\end{equation}

The importance of this equation lies in the fact that the matrix $\mathit{C}$ is only once ​​calculated because it only depends on $x_0$, $h$ and $n$. It does not depend on the function to interpolate. For each data collection $\mathit{Y}$ the only new calculation is the scalar product $\mathit{C} \cdot \mathit{Y}$.

\section{CONVERGENCE ANALYSIS}

Now, the convergence of the developed spline is studied. 

Let $D=\max_{x\in[a,b]} |y(x)-S(x)|$.

Recalling (\ref{1}), results: 
\begin{equation} \label{acot}
D \leq \max_{k \in \{1,\ldots,n\}} \max_{x\in I_k} \left|y-\left(\frac{x-x_{k-1}}{h}y_k - \frac{x-x_k}{h}y_{k-1}\right)\right|+|a_k| h^2.
\end{equation}

From (\ref{ak}) we have:
\begin{equation} \nonumber
|a_k| \leq \sum_{j=1}^{n-1} \left|\frac{j}{n} + s_j - 1\right| |\Delta_j|.
\end{equation}

Noting that $|\frac{j}{n} + s_j - 1| \leq 1 \; \forall j=1,\ldots,n-1$ and assuming that the second derivative of $y(x)$ is bounded, namely $|\Delta_j| \leq M$ for $h$ sufficiently small, results:
\begin{equation} \nonumber
|a_k| \leq (n-1) M = \left(\frac{b-a}{h}-1\right) M.
\end{equation}

In this way, back to (\ref{acot}) we get:
\begin{equation} \nonumber
D \leq \max_{k \in \{1,\ldots,n\}} \max_{x\in I_k} \left|y-\left(\frac{x-x_{k-1}}{h}y_k - \frac{x-x_k}{h}y_{k-1}\right)\right|+ \left(\frac{b-a}{h}-1\right) M h^2.
\end{equation}

Using the result of \cite{kincaid91}:
\begin{equation} \nonumber
\left|y-\left(\frac{x-x_{k-1}}{h}y_k - \frac{x-x_k}{h}y_{k-1}\right)\right| \leq \frac{h^2}{2} \max_{x \in I_k} |y''(x)|.
\end{equation}

It results:
\begin{equation} \nonumber
D \leq M h \left(b-a - \frac{h}{2}\right).
\end{equation}

In this way the convergence order obtained is $o(h)$ for functions with bounded second derivative.

\section{PARITY CONSERVATION}

Next, it will be shown that the developed spline preserves the parity of the interpolated function.

\textbf{Theorem}: Coefficients $a_k$ of the spline $S$ that interpolates $y=y(x)$,
$-a<x<a$, in $n+1$ nodes, verify:

1) If $y$ is even, then $a_k=a_{n-k+1}, \; k=1,\ldots,n$.

2) If $y$ is odd, then $a_k=-a_{n-k+1}, \; k=1,\ldots,n$.

\medskip

\textbf{Demonstration}:

We will do the demonstration for $y$ even.

Since $y$ is even, then $y_j=y_{n-j}, \; j=1,\ldots,n-1$ with what is
immediate that $\Delta_{j}=\Delta_{n-j}, \; j=1,\ldots,n-1$.

\medskip

\underline{Case 1:} $n$ odd.

We will prove by induction on $k$ that $a_k=a_{n-k+1}, \; k=1,\ldots,\frac
{n-1}{2}$.

1) Let's see that $a_1=a_n$. From (\ref{4445}) we have:

\[
a_n=(-1)^{n+1}(a_1+\sum\limits_{j\,=\,1}^{n-1}(-1)^j\Delta_j)=a_1+\sum\limits_{j=1}^{\frac{n-1}{2}}(-1)^j\Delta_j+\sum
\limits_{j=\frac{n+1}{2}}^{n-1}(-1)^j\Delta_{j}=
\]

\[
=a_1+\sum\limits_{j=1}^{\frac{n-1}{2}}(-1)^j\Delta_j+\sum
\limits_{j=\frac{n+1}{2}}^{n-1}(-1)^j\Delta_{n-j}=a_1+\sum
\limits_{j=1}^{\frac{n-1}{2}}(-1)^j\Delta_{j}+\sum\limits_{j=1}
^{\frac{n-1}{2}}(-1)^{n-j}\Delta_j=a_1.
\]

2) Suppose that $a_k=a_{n-k+1}$ for some $k=1,\ldots,\frac{n-3}{2}$.

Let's see that $a_{k+1}=a_{n-k}$. From (\ref{4bis}):
\begin{equation} \nonumber
a_{k+1}=\Delta_k-a_k=\Delta_{n-k}-a_{n-k+1}=a_{n-k}.
\end{equation}

Therefore $a_k=a_{n-k+1}, \; k=1,\ldots,\frac{n-1}{2}$ if $n$ is odd.

\medskip

\underline{Case 2:} $n$ odd.

We will prove by induction on $k$ that $a_k=a_{n-k+1}, \; k=1,\ldots,\frac
{n}{2}$.

1) Let's see that $a_1=a_n$. From (\ref{4445}) we have:
\begin{equation} \nonumber
a_n=(-1)^{n+1}(a_1+\sum\limits_{j=1}^{n-1}(-1)^j\Delta_j
)=-a_1-\sum\limits_{j=1}^{n-1}(-1)^j\Delta_j.
\end{equation}

Then it is enough to prove that $a_1=-\frac{1}{2}\sum\limits_{j=1}
^{n-1}(-1)^j\Delta_j$. From (\ref{666}) we have:

\[
a_1=-\frac{1}{n}\sum\limits_{j=1}^{n-1}(n-j)(-1)^j\Delta
_j=-\frac{1}{2}\{\sum\limits_{j=1}^{n-1}(-1)^j\Delta_j
+\sum\limits_{j=1}^{n-1}[\frac{2(n-j)}{n}-1](-1)^j\Delta_j\}=
\]

\[
=-\frac{1}{2}\{\sum\limits_{j=1}^{n-1}(-1)^j\Delta_j+\sum
\limits_{j=1}^{\frac{n}{2}-1}(1-\frac{2j}{n})(-1)^j\Delta_j
+\sum\limits_{j=\frac{n}{2}+1}^{n-1}(1-\frac{2j}{n})(-1)^j
\Delta_{n-j}\}=
\]

\[
=-\frac{1}{2}\{\sum\limits_{j=1}^{n-1}(-1)^j\Delta_j+\sum
\limits_{j=1}^{\frac{n}{2}-1}(1-\frac{2j}{n})(-1)^j\Delta_j
+\sum\limits_{j=1}^{\frac{n}{2}-1}(-1+\frac{2j}{n})(-1)^{n-j}
\Delta_j\}=-\frac{1}{2}\sum\limits_{j=1}^{n-1}(-1)^j\Delta_j.
\]

2) Suppose that $a_k=a_{n-k+1}$ for some $k=1,\ldots,\frac{n}{2}-1$.

Let's see that $a_{k+1}=a_{n-k}$. From (\ref{4bis}):
\begin{equation} \nonumber
a_{k+1}=\Delta_k-a_k=\Delta_{n-k}-a_{n-k+1}=a_{n-k}.
\end{equation}

Therefore $a_k=a_{n-k+1}, \; k=1,\ldots,\frac{n}{2}$ if $n$ is even.

The proof for $y$ odd and $n$ even is similar to case 1 and the proof for
$y$ odd and $n$ odd is similar to case 2. $\hfill \square$ 

\medskip

\textit{Observation}: In two of the cases ($y$ even, $n$ odd and $y$ odd, $n$ even)
the property is verified independently of the election of
$a_1$, while in the other two cases ($y$ even, $n$ even and $y$ odd,
$n$ odd) the election of $a_1$ is fundamental. It is easy to verify that
for these cases, the property is not verified for every $a_1$.

\medskip

\textbf{Corollary}: The interpolator polynomial $P_k(x)$ of the spline $S$ that interpolates
$y=y(x), -a<x<a$, with $n+1$ nodes, verify:

1) If $y$ is even, then $P_k(x)=P_{n-k+1}(-x), \; k=1,2,\ldots,n$.

2) If $y$ is odd, then $P_k(x)=-P_{n-k+1}(-x), \; k=1,2,\ldots,n$.

\medskip

\textbf{Demonstration}:

Again we will do the demonstration for $y$ even.

Recalling (\ref{4}) is easy to verify that $p_k(x)=p_{n-k+1}(-x) \; k=1,2,\ldots,n$.

Then, taking into account that $a_k=a_{n-k+1}, \; k=1,2,\ldots,n$ and that
$x_k=-x_{n-k}, \; k=0,\ldots,n$, the result is easily obtained.

The demonstration for the case $y$ odd is similar. $\hfill \square$ 

\section{FREDHOLM INTEGRAL EQUATION}

There exist numerous works to determine the numerical solution of Fredholm linear integral equations of the first and second
kind \cite{maleknejad06, maleknejad05, wang14, panda15, chen13}. Following, the results of the previous sections are applied to numerically solve the linear problem. Consider the following equation:

\begin{equation} \label{100}
y(x)-\lambda \int\limits_{a}^{b}K(x,s)y(s) \, ds=f(x),
\end{equation}

\noindent where $y(x)$ is to be determined. When $f(x) \equiv 0$, the equation is homogeneous and $\lambda$ becomes an eigenvalue. Consider also that
$x$, $s$ $\in \mathbb{R}$, $f$ is continuous in $I=[a,b]$, and the kernel $K(x,s)$ is continuous in the region $\Omega=[a,b]\times[a,b]$ with $a,b<\infty$. Let us consider a partition of the interval $I$ with nodes $x_j=a+jh$, $h=\frac{b-a}{n}$, $j=0,1,2,\ldots,n$. Using the quadratic spline $S$ developed in section 2, $y(x)$ is interpolated in $\{y_j\}_{j=0}^{n}$, whose values ​​must be determined. Evaluating (\ref{100}) in nodes $x_{j}$ we get:
\begin{equation} \nonumber
y_j-\lambda{\displaystyle\int\limits_{a}^{b}}K(x_j,s)S(s) \, ds = f_j, \; j=0,1,\ldots,n,
\end{equation}

\noindent where $f_{j}=f(x_{j}).$ Considering (\ref{1}):
\begin{equation} \nonumber
{\displaystyle\int\limits_{a}^{b}}K(x_{j},s)S(s) \, ds={\displaystyle\sum\limits_{k=1}^{n}}\,{\displaystyle\int\limits_{x_{k-1}}^{x_k}}K(x_{j},s)P_{k}(s) \, ds.
\end{equation}

For being $a_k$ linear in $y_j$, $j=0,\ldots,n$, this leads us to a linear system. Taking into account (\ref{4}) and (\ref{1313}) it arises:
\begin{equation} \nonumber
y_j-\lambda{\displaystyle\sum\limits_{k=1}^{n}}\left\{{\displaystyle\sum\limits_{i=0}^{n}}(\delta_{i,k} \, m_{j,k-1}-\delta_{i,k-1} \, m_{j,k}+d_{k,i}) \, y_i\right\} \, = f_j, \; j=0,1,\ldots,n,
\end{equation}

\noindent where $\delta_{i,k}$ is the Kronecker delta, $c_{k,i}$ is that of (\ref{1313}), $m_{j,k}=\displaystyle\frac{1}{h}\int\limits_{x_{k-1}}^{x_k}
K(x_j,s)(s-s_k) \, ds$ and $d_{k,i}=c_{k,i}\int\limits_{x_{k-1}}^{x_k}K(x_j,s)(s-s_{k-1})(s-s_k) \, ds$.

From here:
\begin{equation} \label{105}
{\displaystyle\sum\limits_{i=0}^{n}}(\delta_{j,i}-\lambda\alpha_{j,i})\,y_i = f_j, \; j=0,1,\ldots,n,
\end{equation}

\begin{equation} \nonumber
\alpha_{j,i}={\displaystyle\sum\limits_{k=1}^{n}}\delta_{i,k}\,m_{j,k-1}-\delta_{i,k-1}\,m_{j,k}+d_{k,i}.
\end{equation}

In the non-homogeneous case, (\ref{105}) represent a linear algebraic system
whose solution determines $\{y_j\}_{j=0}^{n}$ and therefore $S(x)$.

\medskip

When $f(x) \equiv 0$ in $[a,b]$, the system is homogeneous and the problem of eigenvalues $\lambda$ ​​and eigenfunctions $y(x)$ is solved in the usual way.

\section{VOLTERRA INTEGRAL EQUATION}

\subsection{Volterra integral equation of the second kind}

Writing the Volterra integral equation of the second kind as:
\begin{equation} \label{201}
y(x)=f(x)+\lambda\int\limits_a^x K(x,s)y(s)\,ds,
\end{equation}
we are interested in determining $y(x)$, $a \leq x \leq b$, where $K(x,s)$
is continuous in $\Omega=[a,b]\times[a,b]$ and $f(x)$ is continuous
in $I=[a,b].$

\bigskip 

For the numerical resolution of (\ref{201}) let us make a
partition of $I$ using nodes $x_j=a+jh$, $h=(b-a)/n$, $j=0,1,2,\ldots,n$. The quadratic spline $S$ interpolates $y(x)$
in $\{y_j\}_{j=0}^n$, whose values ​​must be determined. Then
evaluating (\ref{201}) in the $n+1$ nodes $x_{j}$ follows:
\begin{align}
y_0 & = f_0, \nonumber \\
y_j & = f_j+\lambda\int\limits_{a}^{x_j}K(x_j,s)y(s)\,ds, \; j=1,2,\ldots,n.\nonumber
\end{align}
Taking into account that $y(x)$ is approximated by $S(x)$:
\begin{equation} \nonumber
y_j=f_j+\lambda\sum\limits_{k=1}^j\int\limits_{x_{k-1}}^{x_k}K(x_j,s)P_k(s)\,ds, \; j=1,2,\ldots,n.
\end{equation}
Similarly to section 6, these equations represent a non-homogeneous system of order $n$.
From its resolution $\{y_{j}\}_{j=0}^{n}$ is determined.

\subsection{Volterra integral equation of the first kind}

Writing the Volterra integral equation of the first kind as:
\begin{equation} \label{300}
f(x)=\int\limits_{0}^{x}K(x,s)y(s)\,ds.
\end{equation}

We can derive (\ref{300}) in the particular case in which $K(x,s)$, $\frac{\partial}{\partial x}K(x,s)$, $f(x)$, $f^\prime(x)$ are continuous in $0 \leq x \leq b$, $0 \leq s\leq x$:
\begin{equation} \nonumber
y(x)=\frac{f^\prime(x)}{K(x,x)}-\int\limits_0^x\frac{1}{K(x,x)}\frac{\partial}{\partial x}K(x,s)y(s)\,ds,
\end{equation}

\noindent $K(x,x)$ does not vanish in $0 \leq x \leq b$, which is a second kind equation.

A method is proposed to solve (\ref{300}) which does not require the continuity conditions of $f(x)$ and $K(x,s)$ nor the
restriction on the zeros of $K(x,x)$. In the same way as in the
previous cases, $y(x)$ is determined in nodes $x_k$, $0 \leq x_k \leq b$,
$x_0=0$, $x_n=b$. Taking into account (\ref{300}), $F(x)$ is defined as:
\begin{equation} \label{301}
F(x)=-f(x)+\int\limits_0^x K(x,s)y(s)\,ds,
\end{equation}

\noindent and $F(x_{k})=0$, $k=0,1,\ldots,n$. If $y(x)$ is the solution, $F(x) \equiv 0$.
Taking into account that $y$ is approximated by $S$, (\ref{301}) writes as:
\begin{equation} \nonumber
F(x)=-f(x)+\sum_{k=1}^{n-1}\int\limits_{x_{k-1}}^{x_k}(x-s)^2 P_k(s) \, ds+\int\limits
_{x_{n-1}}^x(x-s)^2 P_{n}(s) \, ds, \; x_{n-1} \leq x \leq x_n.
\end{equation}

To determine $\{y(x_k)\}_{k=0}^n$, we approximate $y(x)$ by the spline $S(x)$
given by (\ref{4}). Evaluating (\ref{301}) in $x_k$, $k=1,\ldots,n$:
\begin{equation} \nonumber
0=-f_k+\sum\limits_{j=1}^k\int\limits_{x_{j-1}}^{x_j}K(x_k,s)P_j(s)\,ds, \; k=1,\ldots,n,
\end{equation}

\noindent and given the fact that $P_j(s)$ depends linearly on $\{y(x_k)\}_{k=0}^n$ a system of 
$n$ non-homogeneous linear equations is obtained which determines
$\{y(x_k)\}_{k=1}^n$ parametrized in $y_0=y(x_0)$, which is determined
by the following ansatz: $y(x_{0})$ is the value that minimizes
$G=\int\limits_{x_{n-1}}^{x_n}F^2(x) \, dx$.

It is easy to see that it is sufficient that $\frac{\partial}{\partial y_0}G=0$. Derivability of $G$ is guaranteed by the linearity of $P_j$ in $y_0$. 

\section{EXAMPLES}

\subsection{Quadratic spline}

Consider $f(x)=\left\vert x\right\vert$, $-1 \leq x \leq 1$, which interpolated
with the Lagrange polynomial $\mathit{L}_n(x)$ with equidistant nodes. It turns out as it is well known, $\lim_{n\rightarrow\infty}\left\Vert \mathit{L}_n(x)-f(x)\right\Vert
_2\neq 0$, $(\left\Vert\cdot\right\Vert$, Lebesgue norm 2). Let $e_n$ be:
\begin{equation} \label{99}
e_n(\cdot)=\int\limits_{x_0}^{x_n}(\cdot-f(x))^{2} \, dx.
\end{equation}

It can be seen that for $n>15$ the Lagrange interpolation $I$ presents sharp
fluctuations. In Table \ref{tab:n1} the error is calculated.

\begin{table}[h!]
 \centering
\begin{tabular}{|c|c|c|c|}
\hline
$n$ & $10$ & $15$ & $20$\\
\hline
$e_n(I)$ & $8.2 \cdot 10^{-2}$ & $2.8 \cdot 10^{-2}$ & $7.2 \cdot 10^{+2}$\\ 
\hline
\end{tabular}\caption{Errors for $f(x)=\left\vert x\right\vert$ with Lagrange interpolation $I$}
\label{tab:n1}
\end{table}

However, using $S$ the undesired
fluctuations are markedly reduced, as seen in Table \ref{tab:n2}.

\begin{table}[h!]
 \centering
\begin{tabular}{|c|c|c|c|c|}
\hline
$n$ & $10$ & $20$ & $50$ & $100$ \\
\hline
$e_n(S)$ & $2.6 \cdot 10^{-3}$ & $6.6 \cdot 10^{-4}$ & 1.0 $10^{-4}$ & $2.6 \cdot 10^{-5}$ \\ 
\hline
\end{tabular}\caption{Errors for $f(x)=\left\vert x\right\vert$ with the $S$ spline}
\label{tab:n2}
\end{table}

\medskip

As a second example, consider $f(x)=\sin2\pi x$, $-1\leq x\leq1$.
Interpolations are made with the Quadratic Spline
end Condition routine ($Q$) from \cite{behforooz88}, with the natural cubic spline of Mathematica 9.0.1.0 software ($M$) and with the $S$ spline from this work. Results are shown in Table \ref{tab:n3}.

\begin{table}[h!]
 \centering
\begin{tabular}{|c|c|c|c|}
\hline
$n$ & $10$ & $50$ & $100$ \\
\hline
$e_n(Q)$ & $2.0 \cdot 10^{-3}$ & $9.8 \cdot 10^{-7}$ & $1.9 \cdot 10^{-8}$\\
\hline
$e_n(M)$ & $1.0\text{ }10^{-3}$ & $8.2\text{ }10^{-11}$ & $1.8\text{ }10^{-13}$ \\
\hline
$e_n(S)$ & $4.0 \cdot 10^{-4}$ & $9.0 \cdot 10^{-9}$ & $1.0 \cdot 10^{-10}$ \\
\hline
\end{tabular}\caption{Errors $e_n$ for $f(x)=\sin2\pi x$}
\label{tab:n3}
\end{table}

Therefore, the spline presented here is a better approximation than the Quadratic Spline end Condition routine.

\subsection{Fredholm equation}

Consider two extra examples from \cite{krasnov82}. The first of them being:
\begin{equation} \label{krasnov1}
y(x)+2\int\limits_0^1 e^{x-t}y(t)\,dt=2x e^x,
\end{equation}

\noindent whose solution is $y(x)=e^x(2x-\frac{2}{3})$ .The solution obtained through
the use of the quadratic spline gives the results shown in Table \ref{tab:n4}, 
where the error $e_n$ is defined in (\ref{99}) and $E_T=\max\left\vert y(x)-y_n(x)\right\vert$
in $x_0\leq x\leq x_n$.

\begin{table}[h!]
 \centering
\begin{tabular}{|c|c|c|}
\hline
$n$ & $5$ & $10$ \\
\hline
$e_n(S)$ & $6.4 10^{-3}$ & $5.1 10^{-4}$ \\
\hline
$E_T$ & $1.6 10^{-3}$ & $1.9 10^{-5}$ \\
\hline
\end{tabular}\caption{Errors $e_n$ and $E_T$ with the $S$ spline for (\ref{krasnov1})}
\label{tab:n4}
\end{table}

The second example corresponds to an equation of the first kind:
\begin{equation} \label{krasnov2}
y(x)-\lambda\int\limits_0^1(2xt-4x^2)y(t)\,dt=0,
\end{equation}

\noindent whose solution is $y(x)=x(1-2x)$, being the eigenvalue $\lambda=-3$, of multiplicity 2. 
It is obatined what is shown in Table \ref{tab:n5}.

\begin{table}[h!]
 \centering
\begin{tabular}{|c|c|c|c|}
\hline
$n$ & $5$ & $9$ & $11$ \\
\hline
$\lambda$ & $-3.21785$ & $-3.065060$ & $-3.04336$ \\
\hline
$e_n(S)$ & $2.5 \cdot 10^{-2}$ & $7.0 \cdot 10^{-3}$ & $4.6 \cdot 10^{-3}$ \\
\hline
$E_T$ & $6.0 \cdot 10^{-2}$ & $1.6 \cdot 10^{-2}$ & $8.6 \cdot 10^{-3}$ \\
\hline
\end{tabular}\caption{Errors $e_n$ and $E_T$ with the $S$ spline for (\ref{krasnov2})}
\label{tab:n5}
\end{table}

Finally, consider from \cite{wang14}:
\begin{equation} \nonumber
\begin{array}{ll}
y(x)  &  =f(x)+\int\limits_0^1(x+s)y(s)\,ds,\\
f(x)  &  =1+\cos(x)-(1+x)\sin1-\cos1,
\end{array}
\end{equation}

\noindent whose solution is $y(x)=\cos x$. $y$ is approximated by $\{1,x,x^2\}$ using the least squares method ($L$) from \cite{wang14} and it is obtained that
$e_2(L)=1.49 \cdot 10^{-3}$. With the $S$ spline and with less computational effort, for $n=5$ and $n=10$, $e_5(S)=1.17 \cdot 10^{-3}$ and $e_{10}(S)=1.79 \cdot 10^{-5}$ are obtained respectively.

\subsection{Volterra equation}

Let us consider from \cite{krasnov82}:
\begin{equation} \label{krasnov3}
y(x)=\frac{1}{1+x^2}-\int\limits_0^x\frac{s}{1+x^2}y(s) \, ds.
\end{equation}

It is obtained what is shown in Table \ref{tab:n6}.

\begin{table}[h!]
 \centering
\begin{tabular}{|c|c|c|}
\hline
$n$ & $5$ & $10$ \\
\hline
$e_n(S)$ & $9.3 10^{-4}$ & $1.8 10^{-4}$ \\
\hline
$E_T$ & $6.1 10^{-5}$ & $5.1 10^{-6}$ \\
\hline
\end{tabular}\caption{Errors $e_n$ and $E_T$ with the $S$ spline for (\ref{krasnov3})}
\label{tab:n6}
\end{table}

From \cite{maleknejad05}, let us consider:
\begin{equation} \label{malek1}
\exp(-x^2)+\frac{\lambda}{2}x(1-\exp(-x^2))v=y(x)+\lambda
\int\limits_0^x xsy(s)\,ds,\;\; 0\leq x\leq 1,
\end{equation}

\noindent for $\lambda=1$, $y(x)=\exp(-x^2)$. The results are shown in Table \ref{tab:n7}.

\begin{table}[h!]
 \centering
\begin{tabular}{|c|c|c|}
\hline
$n$ & $5$ & $10$ \\
\hline
$e_n(S)$ & $4.8 10^{-4}$ & $1.3 10^{-4}$ \\
\hline
$E_T$ & $5.0 10^{-5}$ & $5.1 10^{-6}$ \\
\hline
\end{tabular}\caption{Errors $e_n$ and $E_T$ with the $S$ spline for (\ref{malek1})}
\label{tab:n7}
\end{table}

In the cited work, the error is not specified and the solutions are only compared graphically. For the developed spline, the curves corresponding to the numerical and analytical solutions are completely overlapped.

Again from \cite{krasnov82}:
\begin{equation} \label{krasnov4}
-x^3+\int_0^x(x-s)^2y(s)\,ds=0,
\end{equation}

\noindent whose solution is $y(x)\equiv 3$. It is obtained what is shown in Table \ref{tab:n8}.

\begin{table}[h!]
 \centering
\begin{tabular}{|c|c|c|}
\hline
$n$ & $5$ & $10$ \\
\hline
$e_n(S)$ & $1.4 10^{-9}$ & $4.8 10^{-8}$ \\
\hline
$E_T$ & $6.3 10^{-9}$ & $2.5 10^{-7}$ \\
\hline
\end{tabular}\caption{Errors $e_n$ and $E_T$ with the $S$ spline for (\ref{krasnov4})}
\label{tab:n8}
\end{table}

\section{CONCLUSIONS}

Using variational calculus, a quadratic spline method has been developed that minimizes the spline fluctuations. Spline coefficients are determined explicitly with simple algebraic calculations without needing recursive methods nor solving algebraic systems. The spline has a convergence of order $o(h)$ and when the function has a defined parity, it is conserved by the spline. Given the simplicity of the method, it is useful for solving Volterra-Fredholm integral equations. In a second stage we will extend the results to fractional ordinary differential equations.

\section*{ACKNOWLEDGMENT} This work has been partially sponsored by the Universidad Nacional de Rosario 
through the project ING495 ``Estudio de diversos problemas con ecuaciones diferenciales fraccionarias''. 
The first author is also sponsored by CONICET through an internal doctoral fellowship.


\section{REFERENCES}

\begin{thebibliography}{99}      

\bibitem {kincaid91} D.R. Kincaid, W. Cheney, Numerical Analysis:
Mathematics of Scientific Computing, third ed., Brooks Cole Publishing Company,
California, 1991.

\bibitem {henrici64} P. Henrici, Elements of Numerical Analysis, John Wiley \&
Sons Inc., New York, 1964.

\bibitem {press86} W. Press, S. Teukolsky, W. Vetterling, B. Flannery, Numerical Recipes: The Art of Scientific Computing, first ed.,
Cambridge University Press, New York, 1986.

\bibitem {foucher09} F. Foucher, P. Sablonni\`{e}re, Quadratic spline
quasi-interpolants and collocation methods, Math. Comput. Simulat. 79 (2009) 3455-3465. \\
https://doi.org/10.1016/j.matcom.2009.04.004.

\bibitem {rana89} S.S. Rana, Quadratic Spline Interpolation, J. 
Approx. Theory. 57 (1989) 300-305. \\
https://doi.org/10.1016/0021-9045(89)90045-2.

\bibitem {maleknejad06} K. Maleknejad, H. Derili, Numerical solution of integral
equations by using combination of Spline-collocation method and Lagrange
interpolation, Appl. Math. Comput. 175 (2006) 1235-1244. \\
https://doi.org/10.1016/j.amc.2005.08.034.

\bibitem {maleknejad05} K. Maleknejad, N. Aghazadeh, Numerical solution of Volterra integral equations
of the second kind with convolution kernel by using Taylor-series expansion
method, Appl. Math. Comput. 161 (2005) 915-922. \\
https://doi.org/10.1016/j.amc.2003.12.075.

\bibitem {wang14} Q. Wang, K. Wang, S. Chen, Least squares approximation method for
the solution of Volterra-Fredholm integral equations, J. Comput. 
Appl. Math. 272 (2014) 141-147. \\
https://doi.org/10.1016/j.cam.2014.05.010.

\bibitem {panda15} S. Panda, S.C. Martha, A. Chakrabarti, A modified approach to numerical
solution of Fredholm integral equations of the second kind, Appl. 
Math. Comput. 271 (2015) 102-112. \\
https://doi.org/10.1016/j.amc.2015.08.111.

\bibitem {chen13} F. Chen, P.J.Y. Wong, Solutions of Fredholm integral equations via discrete biquintic splines, Math. 
Comput. Model. 57 (2013) 551-563. \\
https://doi.org/10.1016/j.mcm.2012.07.007.

\bibitem {behforooz88} G. Behforooz, Quadratic Spline, Appl. Math. Lett. 1 (1988) 177-180. \\
https://doi.org/10.1016/0893-9659(88)90067-5.

\bibitem {krasnov82} M.L. Krasnov, A.I. Kiseliov, G.I. Mak\'{a}renko, Integral Equations, Mir, Moscow, 1982.

\end{thebibliography}
\end{document}